\documentclass[10pt]{article}

 \usepackage{amsmath,amsthm,amssymb}
 \usepackage{makeidx,epsfig,lscape}
 \usepackage{color,colortbl}
 \usepackage{fancyhdr}
 \usepackage{xcolor,pict2e}
\usepackage{graphicx} %package to manage images

 \renewcommand{\headrulewidth}{0pt}
 \renewcommand{\footrulewidth}{0.5pt}

 \definecolor{myaqua}{rgb}{0.0,0.5,0.55}
 \definecolor{lightaqua}{rgb}{0.75,0.95,0.95}

 \usepackage[colorlinks = true,
            linkcolor = myaqua,
            urlcolor  = blue,
            citecolor = myaqua]{hyperref}

\usepackage{caption}
\usepackage{floatrow}
\def\lin#1#2{\textcolor[rgb]{0.6,0.6,0.6}{\vspace*{#1mm} \hrule
   height 3 pt \vspace*{#2mm}}}
\def\bt{\begin{tabular}}
\def\et{\end{tabular}}
\def\and{\mbox{ and }}

\def\1{{\bf 1}}

 \def\sectionn#1{\refstepcounter{section}{\color{myaqua}

 \vskip 6mm

 \noindent\Large\bf\thesection. #1}

 \vskip 3mm}

 \def\boxx#1#2#3#4#5{
 {\linethickness{#4pt}\put(#1,#5){\color{myaqua}{\line(1,0){#3}}}}
 \multiput(#1,#2)(0,#4){2}{\line(1,0){#3}}
 \multiput(#1,#2)(#3,0){2}{\line(0,1){#4}}
  }

\begin{document}

 \fancyhead[L]{\hspace*{-13mm}
 \bt{l}{\bf Open Journal of *****, 2014, *,**}\\
 Published Online **** 2014 in SciRes.
 \href{http://www.scirp.org/journal/*****}{\color{blue}{\underline{\smash{http://www.scirp.org/journal/****}}}} \\
 \href{http://dx.doi.org/10.4236/****.2014.*****}{\color{blue}{\underline{\smash{http://dx.doi.org/10.4236/****.2014.*****}}}} \\
 \et}
 \fancyhead[R]{\includegraphics{pic1.ps}}

 $\mbox{ }$

 \vskip 12mm

{ % \fontfamily{Cambria}\selectfont

% "Title of the Paper"
{\noindent{\huge\bf\color{myaqua}
Obtaining a New Representation for the Golden\\[2mm]  
Ratio by Solving a Biquadratic Equation}}
%
% \runtitle{Obtaining a New Representation for the Golden Ratio by Solving a Biquadratic Equation}
\\[6mm]
{\large\bf Leonardo Mondaini$^{1,2}$}}
\\[2mm]
{ %\fontfamily{Calibri}\selectfont
 $^1$Department of Oncology, University of Alberta, Edmonton AB, Canada\\$^2$Grupo de F\'{i}sica Te\'orica e Experimental, Departamento de Ci\^encias Naturais,\\Universidade Federal do Estado do Rio de Janeiro, Rio de Janeiro RJ, Brazil\\
Email:
\href{mailto:mondaini@ualberta.ca}{\color{blue}{\underline{\smash{mondaini@ualberta.ca}}}},
\href{mailto:mondaini@unirio.br}{\color{blue}{\underline{\smash{mondaini@unirio.br}}}}
 \\[4mm]

\lin{5}{7}

 { % \fontfamily{Cambria}\selectfont
 {\noindent{\large\bf\color{myaqua} Abstract}{\bf \\[3mm]
 \textup{In the present work we show how different ways to solve biquadratic
equations can lead us to different representations of its solutions.
A particular equation which has the golden ratio and its reciprocal
as solutions is shown as an example.
 }}}
 \\[4mm]
 {\noindent{\large\bf\color{myaqua} Keywords}{\bf \\[3mm]
 Golden Ratio; Algebraic Equations; Recreational Mathematics; History of Mathematics
}

 \fancyfoot[L]{{\noindent{\color{myaqua}{\bf How to cite this
 paper:}} L. Mondaini (2014)
  Obtaining a New Representation for the Golden Ratio by Solving a Biquadratic Equation.
 ***********,*,***-***}}

\lin{3}{1}

\sectionn{Introduction}

{ \fontfamily{times}\selectfont
 \noindent 
The study of algebraic equations has occupied the brightest mathematical minds
throughout many centuries. We must highlight among the main results of the studies in
this area, the attainment of the formula for resolution of the general quadratic equations which, along with the formula for resolution of the general cubic equations achieved by Niccolo Fontana (Tartaglia)\footnote{Usually the formula for resolution of the general cubic equation is attributed to Girolamo
Cardano, and being thus it receives his name. A more detailed description about the dispute for
the priority of this solution can be found in \cite{boyer}.}, led to the creation of complex numbers,
since the application of these formulas led to a ``misterious" thing: the square root
of a negative number.
The solution of general quartic equations by Ludovico Ferrari (a pupil of Cardano) comes to complete this picture, once it was established that a solution by
radicals for generic equations of degree equal to or greater than 5 cannot be achieved
(a result proved for the first time by the prodigies Niels Abel and Evariste Galois). In the present work, which may be classified into the field of recreational mathematics and is devoted to stimulate the interest of readers with pre-university level mathematical background as a way of inspiring their further study on the subject, we will focus on an interesting aspect associated to a particular kind of quartic equation, namely, the biquadratic one.

\renewcommand{\headrulewidth}{0.5pt}
\renewcommand{\footrulewidth}{0pt}

 \pagestyle{fancy}
 \fancyfoot{}
 \fancyhead{} % clear all header and footer fields
 \fancyhf{}
 \fancyhead[RO]{\leavevmode \put(-90,0){\color{myaqua}L. Mondaini} \boxx{15}{-10}{10}{50}{15} }
 \fancyhead[LE]{\leavevmode \put(0,0){\color{myaqua}L. Mondaini}  \boxx{-45}{-10}{10}{50}{15} }
 \fancyfoot[C]{\leavevmode
 \put(0,0){\color{lightaqua}\circle*{34}}
 \put(0,0){\color{myaqua}\circle{34}}
 \put(-2.5,-3){\color{myaqua}\thepage}}

 \renewcommand{\headrule}{\hbox to\headwidth{\color{myaqua}\leaders\hrule height \headrulewidth\hfill}}

The rest of this work is organized as follows. A particular equation which has the golden ratio and its reciprocal
as solutions is presented in Section 2. In Sections 3 and 4, we solve this equation by using two different algorithms. Finally, in Section 5, we present our concluding remarks.

{ \fontfamily{times}\selectfont
 \noindent

\sectionn{An interesting equation}
\label{sec:K-M}

{ \fontfamily{times}\selectfont
 \noindent
We start by considering the following characteristic equation \cite{wolfram,seymour}
\begin{equation}
\rm{det}(\mathcal{A}-\lambda\mathcal{I})=0, \label{eq1}
\end{equation}
where $\mathcal{A}$ is a $4 \times 4$ symmetric real matrix defined by
\begin{equation}
\mathcal{A}\equiv \begin{pmatrix}
0& 1& 0& 0\\
1& 0& 1& 0\\
0& 1& 0& 1\\
0& 0& 1& 0
\end{pmatrix}, \label{eq2}
\end{equation}
and $\mathcal{I}$ is the corresponding $4 \times 4$ identity matrix. Notice that we can rewrite equation (\ref{eq1}) in the following way
\begin{equation} 
\begin{vmatrix}
-\lambda& 1& 0& 0\\
1& -\lambda& 1& 0\\
0& 1& -\lambda& 1\\
0& 0& 1& -\lambda
\end{vmatrix}=0, \label{eq3}
\end{equation}
or
\begin{equation}
\lambda^4-3\lambda^2+1=0. \label{eq4}
\end{equation}
This simple biquadratic equation displays an interesting feature. The form its
solutions are expressed depends, apparently, on the algorithm used for solving it.
Even more interesting is the fact that one of these algorithms yields the numbers
$(\sqrt{5}+1)/2$ and $(\sqrt{5}-1)/2$ (as well as their respective symmetrical ones) as solutions, which
are the most known representation of this ubiquitous mathematical phenomenon, namely, the {\it golden ratio} \cite{livio}, and its reciprocal. This famous number appears historically as the
solution of the quadratic equation 
\begin{equation}
x(x+a)=a^2,\label{eq5}
\end{equation}
which is related to the geometrical problem of dividing a given line segment $\overline{AB}$
into what is called the {\it golden section}, which the celebrated astronomer Johannes
Kepler called ``one of the two Jewels of Geometry" (the second one being the Pythagorean theorem). Translated into mathematical language, the golden section means that
the segment $\overline{AB}=a$ is cut at a point $C$ so that the whole segment is in the same ratio to the larger part $\overline{AC }= x$ as $\overline{AC }$ is to the other part, $\overline{CB}= a-x$. That is
\begin{equation}
\frac{a}{x}=\frac{x}{a-x};\,\,\,\,\,\,\,\,\,x>a-x. \label{eq6}
\end{equation}
This, in turn, leads to the quadratic equation $x(x + a) = a^2$ already mentioned,
the positive root of which is $x=a(\sqrt{5}-1)/2$. Notice that when $a=1$, the value $x=(\sqrt{5}-1)/2$
is the reciprocal of the golden ratio, i.e., 0.6180339\dotso
In the next sections we will apply two different algorithms in order to solve the equation (\ref{eq4}).
{ \fontfamily{times}\selectfont
 \noindent

\sectionn{First Algorithm}
\label{sec:K-M2}

{ \fontfamily{times}\selectfont
 \noindent
Firstly, we will solve equation (\ref{eq4}) by using an algorithm very similar to the one originally employed
by Ludovico Ferrari in his solution for the quartic equations (polynomial equations of
the fourth degree)\footnote{The goal of Ferrari's algorithm for solving the general quartic is to have perfect squares in
both sides of the equation.}\cite{cardano}. In order to do that, we start by observing that equation (\ref{eq4}) can be also
rewritten in the following way:
\begin{equation}
\lambda^4-3\lambda^2+1=\lambda^4-2\lambda^2-\lambda^2+1=0,
\label{eq7}
\end{equation}
or
\begin{equation}
\lambda^4-2\lambda^2+1 = \lambda^2.
\label{eq8}
\end{equation}
The left side of this equation, a perfect square, may be trivially simplified as 
\begin{equation}
(\lambda^2-1)^2 = \lambda^2,
\label{eq9}
\end{equation}
which implies that
\begin{equation}
\lambda^2-1 = \pm \lambda,
\label{eq10}
\end{equation}
or
\begin{equation}
\lambda^2 \mp \lambda-1= 0.
\label{eq11}
\end{equation}
When solving the above equations by using the well-known quadratic formula, we find that the solutions for the equation with $-\lambda$ are given by
\begin{equation}
\lambda_-=\begin{cases}
\frac{\sqrt{5}+1}{2}\\
-\frac{\sqrt{5}-1}{2}
\end{cases},
\label{eq12}
\end{equation}
whereas for the equation with $+\lambda$ (which is identical to equation (\ref{eq5}) with $a=1$) we have the following solutions
\begin{equation}
\lambda_+=\begin{cases}
\frac{\sqrt{5}-1}{2}\\
-\frac{\sqrt{5}+1}{2}
\end{cases}.
\label{eq13}
\end{equation}

Thus, the complete set of solutions of the original biquadratic equation is given by
\begin{equation}
\lambda=\begin{cases}
\frac{\sqrt{5}+1}{2}\\
-\frac{\sqrt{5}+1}{2}\\
\frac{\sqrt{5}-1}{2}\\
-\frac{\sqrt{5}-1}{2}
\end{cases},
\label{eq14}
\end{equation}
where we remind again the reader that $\frac{\sqrt{5}+1}{2}$ is the usual representation for the golden
ratio, and $\frac{\sqrt{5}-1}{2}$ is its reciprocal.

{ \fontfamily{times}\selectfont
 \noindent

\sectionn{Second Algorithm}
\label{sec:K-M3}

{ \fontfamily{times}\selectfont
 \noindent
The second algorithm follows the conventional method to solve biquadratic equations. In such method we start by considering the following change of variables
\begin{equation}
\eta=\lambda^2
\label{eq15}
\end{equation}
which allows us to reduce the equation (\ref{eq4}) to the form
\begin{equation}
\eta^2-3\eta+1=0.
\label{eq16}
\end{equation}
A direct application of the quadratic formula yields the solutions
\begin{equation}
\eta=\begin{cases}
\frac{\sqrt{5}+3}{2}\\
-\frac{\sqrt{5}-3}{2}
\end{cases}.
\label{eq17}
\end{equation}
Since $\lambda=\pm\sqrt{\eta}$, we then have the following set of solutions for the original equation\footnote{The solutions are displayed so as to indicate a direct correspondence with the set of solutions
obtained with the first algorithm.}
\begin{equation}
\lambda=\begin{cases}
\sqrt{\frac{\sqrt{5}+3}{2}}\\
-\sqrt{\frac{\sqrt{5}+3}{2}}\\
\sqrt{\frac{3-\sqrt{5}}{2}}\\
-\sqrt{\frac{3-\sqrt{5}}{2}}
\end{cases},
\label{eq14}
\end{equation}

\sectionn{Concluding Remarks}

{ \fontfamily{times}\selectfont
 \noindent
We saw that apparently distinct solutions are obtained by solving equation (\ref{eq4}) by
two different algorithms. However, we can easily verify that they
are identical. Just compute the square of $(\sqrt{5}+1)/2$ and check it!
We have then obtained a new representation for the golden ratio, namely
\begin{equation}
\rho=\frac{\sqrt{5}+1}{2}=\sqrt{\frac{\sqrt{5}+3}{2}}.
\label{eq14}
\end{equation}

 {\color{myaqua}

 \vskip 6mm

 \noindent\Large\bf Acknowledgments}

 \vskip 3mm

{ \fontfamily{times}\selectfont
 \noindent
This work has been supported in part by CNPq.

{\color{myaqua}

}}
\end{document}